\documentclass[12pt]{article}
\usepackage{amsfonts}
\usepackage{latexsym}
\usepackage{graphics}
\begin{document}

\title{Restrictions of Harmonic Functions on the Sierpinski Gasket to Segments
\thanks{Math Subject Classifications. Primary 31C45, 42C99.}}

\author{B\"{u}nyamin Demir, Vakif Dzhafarov,\\ \c{S}ahin Ko\c{c}ak, Mehmet \"{U}reyen}

\date{}
\maketitle
\begin{abstract}
The restrictions of a harmonic function on the Sierpinski
Gasket (SG) to the segments in SG have been of some interest. We show
that the sufficient conditions for the monotonicity of these restrictions
given by Dalrymple, Strichartz and Vinson are also necessary. We then
prove that the normal derivative of a harmonic function on SG on the
junction points of the contour of a triangle in SG is always nonzero
with at most a single exception.

We finally give an explicit derivative computation for the restriction
of a harmonic function on SG to segments at specific points of the
segments: The derivative is zero at points dividing the segment in ratio 1:3.
This shows that the restriction of a harmonic function to a segment
of SG has the following curious property: The restriction has infinite
derivatives on a dense set of the segment (at junction points) and
vanishing derivatives on another dense set.
\end{abstract}

\centerline{\bf 1. Introduction}

We will first briefly recall the rudiments of harmonic analysis on the
Sierpinski Gasket. ([1],[2],[3],[4])

Let $K$ be the Sierpinski Gasket (SG) constructed on the unit equilateral
triangle $G_0$ with vertices $\{ p_0, \, p_1, \, p_2 \, \}$ and
$G_m$ be the graph in the $m$ th step as in the following figure.\\

\begin{figure}
\begin{center}
\includegraphics{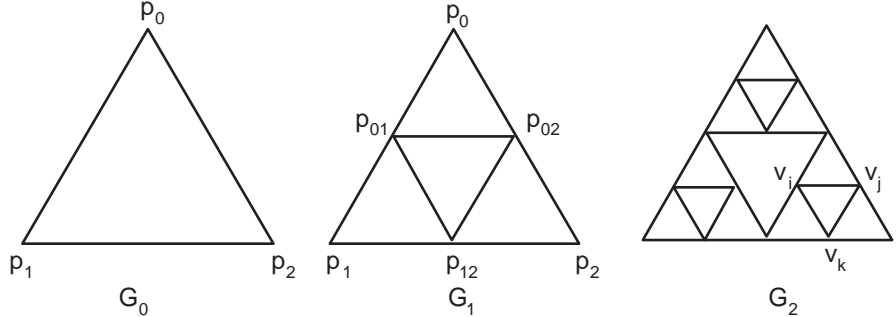}
\caption{Iterated graphs in SG}
\end{center}
\end{figure}

\noindent
{\bf Definition.}
{\it The function $f \in C(K), \, f: K \longrightarrow \mathbb{R}$
is called harmonic on $K$ if for every minimal triangle in $G_m \, (m\geq 1)$,
with vertices $ \{ v_i, \, v_j, \, v_k \}$ the equalities
\begin{equation}
f(v_i)+f(v_j)+f(v_{ik})+f(v_{jk})-4f(v_{ij})=0
\end{equation}
\noindent
hold, where $v_{ij}$ is the midpoint of the segment $[v_i, v_j]$.}\\

Let
\begin{equation}
f(p_0)=\alpha, \, f(p_1)=\beta, \, f(p_2)=\gamma \, .
\end{equation}

\noindent
Then this triple
$(\alpha, \beta, \gamma)$ completely defines a harmonic function $f$,
that is, there exists a unique harmonic function $f: K \longrightarrow \mathbb{R}$
such that
$f(p_0)=\alpha, \, f(p_1)=\beta$ and  $f(p_2)=\gamma$. This harmonic function
depends linearly on the triple $(\alpha, \beta, \gamma)$. According to
the harmonic extension algorithm, it holds
\begin{equation}
f(p_{12})=\frac{1}{5}(\alpha+2\beta+2\gamma), f(p_{02})=\frac{1}{5}(2\alpha+\beta+2\gamma),
f(p_{01})=\frac{1}{5}(2\alpha+2\beta+\gamma) .
\end{equation}

\noindent
We are interested in restrictions of the harmonic function $f$ to the line
segments contained in the SG. From (1)-(3) it can be seen that, if a nonconstant
harmonic function is monotone on some line segment that is contained in SG,
then it is strictly monotone on it.\\

\noindent
Let $T_m $ be a minimal triangle with vertices $v_i, v_j$ and $v_k$
in $G_m$. The
sides of $T_m$ can be ordered by the values $|f(v_i)-f(v_j)|.$\\

\noindent
{\bf Theorem 1 ([DSV]).}\\
{\it Let $f$ be a harmonic function on SG, E an edge in
$G_m$ with endpoints $v_0, \, v_1$ and midpoint $v_{01}$. Suppose

\begin{equation}
f(v_0)<f(v_{01})<f(v_1)
\end{equation}

\noindent
and

\begin{equation}
\frac{1}{4} \leq \frac{f(v_1)-f(v_{01})}{f(v_{01})-f(v_0)} \leq 4 \, .
\end{equation}

\noindent
Then the restriction of $f$ to E is strictly increasing.}\\

\noindent
{\bf Theorem 2 ([DSV]).}\\
{\it The restriction of $f$ to the two largest edges of $T_m$
is monotone. On the smallest edge of $T_m$, the restriction of $f$ might be
monotone or not; but if it is not monotone, then it has a unique extremum.}\\
(We changed the wording of the Theorem 2 in [DSV] slightly.)\\

In this paper we will show that the sufficient conditions in Theorem 1 are also
necessary for the monotonicity of the restriction to an adge in $G_m$. Then
we  will characterize when exactly simultaneous monotonicity of restrictions
to all three edges of a triangle $T_m$ in $G_m$ occurs. Furthermore, we will
prove that at the junction points of any segment E in $G_m$, the derivatives of the restriction
of a harmonic function exist improperly (possibly with exception at a single point),
and we will also prove that on another dense subset of E the derivatives of the restriction
are zero.\\

\noindent
Now we remark that it is enough to prove these statements for the triangle
$G_0$ instead of considering an arbitrary triangle $T_m$ in $G_m$, because
the procedure of harmonic extension is the same for $G_0$ or $T_m$. In this
spirit and for simplicity we now reformulate the above two theorems of DSV
(Dalrymple, Strichartz and Vinson):\\

\noindent
{\bf Theorem $1^{\prime}$.}\\
{\it Let $f$ be a harmonic function on SG and assume
$\beta = f(p_1)<\gamma = f(p_2)$. Let $p_{12}$ denote the midpoint of
$[p_1, p_2]$. Then, the restriction of $f$ to the edge $[p_1, p_2]$ is
strictly increasing, if the inequalities\\

\noindent
i) $\beta < f(p_{12})<\gamma$\\

\noindent
ii) $\frac{1}{4}\leq \frac{\gamma-f(p_{12})}{f(p_{12})-\beta} \leq 4$\\

\noindent
are satisfied. (We remark that by harmonicity $f(p_{12})=\frac{2\beta+2\gamma+\alpha}{5}$)}\\

\noindent
{\bf Theorem $2^{\prime}$.}\\
{\it Order the edges of $G_0$ by the values $|f(p_0)-f(p_1)|=|\alpha - \beta |, \,
|f(p_0)-f(p_2)|=|\alpha - \gamma |, \, |f(p_1)-f(p_2)|= |\beta - \gamma | $.
Then the restrictions of $f$ to the two largest of $[p_0, p_1], \, [p_0, p_2]$ and $[p_1, p_2]$
are monotone. On the smallest edge, the restriction of $f$ might be monotone
or not; but if it is not monotone, then it has a unique extremum.}\\

Before proceeding further, we will recast the inequalities i) and ii) of
Theorem $1^{\prime}$ in a more sympatetic form:\\

\noindent
{\bf Lemma 1.}\\
{\it The DSV inequalities\\

\noindent
i) $\beta < f(p_{12})<\gamma$\\

\noindent
ii) $\frac{1}{4}\leq \frac{\gamma-f(p_{12})}{f(p_{12})-\beta} \leq 4$\\

\noindent
are equivalent to the inequalities\\

\noindent
iii) $\beta < \gamma, \quad 2\beta - \gamma \leq \alpha \leq 2 \gamma - \beta \, .$}\\

\noindent
(Proof is  straightforward. )\\

\noindent
\centerline{\bf 2. Characterization of the Monotonicity of the Restrictions}\\

Consider the side $[p_1, \, p_2]=[0, \, 1]$ of $G_0$ and the restriction of the
harmonic function $f$ defined by (2) to $[p_1, \, p_2]$. The following
lemma can be proved by induction on $m$.\\

\noindent
{\bf Lemma 2.}\\
{\it Let $l_m=\frac{1}{2}-\frac{1}{2^{m+1}}, \,
r_m=\frac{1}{2}+\frac{1}{2^{m+1}} \quad (m=1, 2, 3, ...).$  Then

\begin{equation}
f(\frac{1}{2^m}) = \frac{3^m-1}{2\cdot5^m}\alpha + [1-(\frac{3}{5})^m]\beta +
\frac{3^m+1}{2\cdot5^m}\gamma
\end{equation}

\begin{equation}
f(1- \frac{1}{2^m}) = \frac{3^m-1}{2\cdot5^m}\alpha + [1-(\frac{3}{5})^m]\gamma +
\frac{3^m+1}{2\cdot5^m}\beta
\end{equation}

\begin{equation}
f(l_m) = \frac{5^m-1}{5^{m+1}}\alpha + \frac{3^{m+1}+4\cdot5^m+3}{10\cdot5^m} \beta +
\frac{4\cdot5^m-3^{m+1}-1}{10\cdot5^m}\gamma
\end{equation}

\begin{equation}
f(r_m) = \frac{5^m-1}{5^{m+1}}\alpha + \frac{3^{m+1}+4\cdot5^m+3}{10\cdot5^m} \gamma +
\frac{4\cdot5^m-3\cdot3^{m+1}-1}{10\cdot5^m}\gamma
\end{equation}
}
\noindent
(Actually, by symmetry, one of these equalities implies the other three).\\

\noindent
{\bf Theorem 3.}\\
{\it Let $f$ be the harmonic function on SG generated by the
triple $ (\alpha, \beta, \gamma )$. Then the restriction of $f$ to
$[p_1, p_2]$ is strictly increasing if and only if $\beta < \gamma$ and
$2 \beta - \gamma \leq \alpha \leq 2 \gamma - \beta$.}\\

\noindent
{\bf Remark 1.} This fact can also be expressed as follows: the restriction
of a nonconstant $f$ to $[p_1, p_2]$ is strictly monotone iff $\alpha$
lies between $2\beta - \gamma$ and $2 \gamma - \beta $, because for a
nonconstant $f$, necessarily $\beta \neq \gamma$ : if we had $\beta = \gamma $,
then $\alpha $ lying between $2 \beta - \gamma $ and $2 \gamma - \beta$
would coincide with $\beta$ and $\gamma$ making the function constant.

\noindent
{\bf Remark 2.} The condition "$\alpha $ lies between $2 \beta - \gamma $ and
$2 \gamma - \beta$" can also be expressed as follows: Let $\delta = \alpha + \beta + \gamma$.
Then $2 \beta - \gamma - \alpha = 3 \beta - \delta , \, 2 \gamma - \beta - \alpha = 3 \gamma - \delta$
and the condition takes the form $(3\beta - \delta )(3 \gamma - \delta ) \leq 0$.\\

\noindent
{\bf Proof of Theorem 3.} If $\beta < \gamma $ and $2\beta - \gamma \leq \alpha \leq 2 \gamma - \beta$,
then by Lemma 1, the DSV inequalities are satisfied, and consequently the restriction
of $f$ to $[p_1, p_2]$ is strictly increasing by Theorem $1^{\prime}$. Now
we show the necessity of the inequalities $\beta < \gamma $ and
$2\beta - \gamma \leq \alpha \leq 2 \gamma - \beta$. The inequality $\beta < \gamma $ is
obvious. For all $m=1,2,...$ we have by Lemma 2 and the inequality assumption\\

\begin{equation}
f(\frac{1}{2^m}) = \frac{3^m-1}{2\cdot5^m}\alpha + [1-(\frac{3}{5})^m]\beta +
\frac{3^m+1}{2\cdot5^m}\gamma > \beta
\end{equation}

\begin{equation}
f(1- \frac{1}{2^m}) = \frac{3^m-1}{2\cdot5^m}\alpha + [1-(\frac{3}{5})^m]\gamma +
\frac{3^m+1}{2\cdot5^m}\beta < \gamma
\end{equation}

\noindent
From (10) and (11) we obtain

\begin{equation}
\sup_m \frac{2\cdot3^m \beta - 3^m \gamma - \gamma }{3^m - 1} \leq
\alpha \leq \inf_m \frac{2\cdot3^m \gamma - 3^m \beta - \beta }{3^m - 1}
\end{equation}

\noindent
Computing the supremum and infimum gives
$$
2 \beta - \gamma \leq \alpha \leq 2 \gamma - \beta \quad \Box
$$

Now we will find conditions on $\alpha, \, \beta$ and $\gamma$ which guarantee
the monotonicity of restrictions of a harmonic function $f$ on SG to all three
edges of $G_0$.\\

\noindent
{\bf Theorem 4.}\\
{\it The restrictions of a nonconstant harmonic function $f$ on SG to all three
edges of the triangle $G_0$ are simultaneously strictly monotone if and
only if one of the following conditions  holds:
$$
2\alpha = \beta + \gamma
$$

\noindent
or

$$
2\beta = \alpha + \gamma
$$

\noindent
or

$$
2\gamma = \alpha + \beta
$$
}
\noindent
{\bf Proof.} By Theorem 3 and the Remarks 1  and 2, the restriction
of $f$ to all three edges of SG are simultaneously strictly monotone if and
only if all of the following three inequalities are satisfied:\\

$$
(3\beta - \delta )(3\gamma - \delta ) \leq 0
$$

$$
(3\gamma - \delta )(3\alpha - \delta ) \leq 0
$$

$$
(3\alpha - \delta )(3\beta - \delta ) \leq 0.
$$

\noindent
These three inequalities can only be satisfied if at least one of the numbers
$3\alpha - \delta , 3\beta - \delta$ or $ 3\gamma - \delta$ vanishes. This is
also sufficient as easily can be seen. For example, let $3\alpha - \delta=0$.
Then the second and third inequalities are obviously satisfied and the first
is also true, because then $2\alpha = \beta + \gamma$ and this means that
$\alpha= \frac{1}{2}[(2\beta - \gamma ) + (2 \gamma - \beta )]$ and this is
enough for the first inequality to be satisfied $\Box$\\

\noindent
{\bf Remark 3.} Interestingly, the equations $2\alpha = \beta + \gamma, \,
2\beta = \alpha + \gamma , \, 2\gamma = \alpha + \beta $
are equivalent to the vanishing of the normal
derivatives at the vertices. [BST]\\

\noindent
{\bf Remark 4.} We can express the content of Theorem 4 also in terms of
side lengths of $G_0$ (with respect to $f$) : By definition, the lengths
of the edges of SG are $|\alpha - \beta |, \, |\alpha - \gamma |, \,
|\beta - \gamma|$. By Theorem 2 ([DSV]) we know that the restriction of $f$
to the two largest edges are monotone. After having characterized
simultaneous strictly monotonicity on all three edges in Theorem 4, we can
say that this is exactly the case if and only if the $G_0$ is isosceles,
the third side being the largest: For example, if $2\alpha = \beta + \gamma$,
then $|\alpha - \beta | = |\alpha - \gamma |$ and $|\beta - \gamma | = 2 |\alpha - \beta | >
|\alpha - \beta | $. In other words, a non-monotone restriction occurs, iff
the length of the smallest edge is strictly less than the lengths of the other
two edges. So, there are two cases  for non-monotonicity: Either $G_0$ is
scalene, or , if it is isosceles, then the length of the third side is strictly
less than the others. In the latter case, the third side of the isosceles
must have necessarily side length zero. (See Fig. 2)\\

\begin{figure}[h]
\begin{center}
\includegraphics{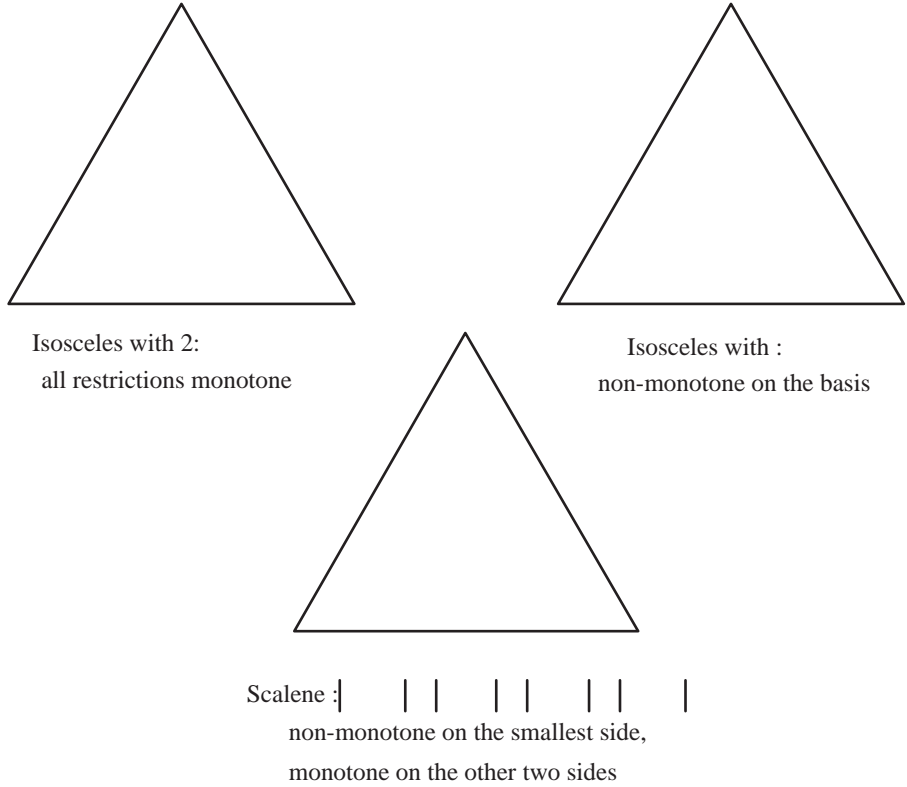}
\caption{Monotonicity classification with respect to side-lengths of $G_0$ }
\end{center}
\end{figure}

\centerline{\bf 3. Derivatives of Restrictions at the Junction
Points}

We will show in this section that the derivatives of the
restriction of a harmonic function on SG to any segment in SG
exist and are infinite at all junction points with possibly a
single exception. It is again enough to show this for the edges of
$G_0$ by general reasons we indicated in the
introduction.\\

We first need the following \\

\noindent
{\bf Lemma 3.}\\
{\it Let the function $g:[0, 1] \longrightarrow \mathbb{R} $
be strictly monotone in a neighborhood of $x_0 \in [0, 1],$  $d \in (0, 1), \, a \neq 0$
and $x_m=x_0+ad^m$. Assume

$$
\frac{g(x_m)-g(x_0)}{x_m-x_0}
$$

\noindent
is defined and tends to $0$ (or $\pm \infty $) as  $m \longrightarrow \infty$.
If $a<0$ then the left derivative of $g$ at $x_0$ exists and is $0$ (or $\pm \infty $);
If $a>0$ then the right derivative of $g$ at $x_0$ exists and is $0$ (or $\pm \infty $).}\\

\noindent
{\bf Proof.} We consider only the case, where $g$ is monotone increasing and
$a>0$. Let $x_0 \in [0, 1)$ and $x>x_0$. Then there exists $m \in \mathbb{N}$
such that

$$
x_0+ad^{m+1} \leq x \leq x_0+ad^m
$$

\noindent
As $x$ tends to $x_0$ , $m$  tends to infinity and from the inequalities

$$
d. \frac{g(x_{m+1})-g(x_0)}{x_{m+1}-x_0} \leq \frac{g(x)-g(x_0)}{x-x_0}
\leq \frac{1}{d}. \frac{g(x_m)-g(x_0)}{x_m-x_0}
$$

\noindent
we get the result $\Box$\\

\noindent
{\bf Remark 5.}  In the above lemma, one-sided monotonicity is obviously
enough for one-sided derivative calculations.\\

\noindent
Now we will compute the derivative of the restriction at the point $p=1/2$, for monotone
restrictions.\\

\noindent
{\bf Lemma 4.}\\
{\it Let the restriction of the harmonic function $f$ to the edge
$[p_1,p_2]=[0,1]$ be strictly monotone. Then $f^{\prime}(\frac{1}{2})=+\infty$
for $f$ monotone increasing and $f^{\prime}(\frac{1}{2})=-\infty$ for $f$
monotone decreasing.}\\

\noindent
{\bf Proof.} We give the proof for $f$ monotone increasing:\\
Applying (8) and Lemma 1, where $x_0=\frac{1}{2}, \, a=-\frac{1}{2}, \,
d=\frac{1}{2}$, we obtain
$$
\lim_{m \longrightarrow \infty} \frac{f(l_m)-f(\frac{1}{2})}{l_m-\frac{1}{2}}
=\lim_{m \longrightarrow \infty} (\frac{3}{5}).(\frac{6}{5})^m(\gamma - \beta )
= + \infty.
$$

\noindent
Then by Lemma 3  the left-hand derivative at $p=\frac{1}{2}$ is $+\infty$.
Analogously, using (9), we obtain that the right-hand derivative at
$p=\frac{1}{2}$ is also $+\infty$ $\Box$\\

Applying Lemma 4 to smaller triangles, we see that the derivatives exist
improperly at all inner junction points of $[p_1,p_2]$ in whose vicinity
the restriction is strictly monotone. Using Lemma 3 and Lemma 2, (6) and (7),
we can compute the derivatives at $p_1$ and $p_2$ also.\\

\noindent
{\bf Lemma 5.}\\
{\it If $2\beta = \alpha + \gamma ,$ then the derivative of the
restriction of $f$ to $[p_1,p_2]$ at $p_1$ is zero; otherwise infinite.
Similarly, if $2\gamma=\alpha + \beta ,$ then the derivative of the
restriction at $p_2$ is zero, otherwise infinite.}\\
(This result is implicit in Theorem 4 of [DSV]).\\

Lemma 4 and 5 yield to \\

\noindent
{\bf Lemma 6.}\\
{\it Let $f$ be a non-constant harmonic function on SG. Then the
derivatives of the restriction of $f$ to $[p_1,p_2]$ exist improperly at
all junction points on $[p_1,p_2]$, with possibly a single exception.
Moreover, this is true for the whole contour $G_0=[p_0,p_1] \cup [p_0,p_2] \cup
[p_1,p_2]$, i.e. the derivative of the restrictions can vanish only for a
single junction point on the whole contour and is infinite for all the junction
points.}\\

The derivatives of the restriction are related to normal derivatives (For normal
derivatives see [BST]). It can be seen from the work of [BST] using monotonicity
of restrictions of a harmonic function $f$ on SG, that if the normal derivative
at a junction point vanishes, then the derivative of the restriction
to a segment containing the junction point also vanishes at that junction
point. But we have seen above, that at junction points on the contour $G_0$
the derivative is infinite with possibly a single exception. This proves
that the normal derivative is nonzero for all junction points on the contour $G_0$
with at most a single exception for a non-constant $f$.\\

\noindent
{\bf Remark 6.} It can be shown that, if the numbers $\alpha, \beta, \gamma$
(being not all equal) are linearly independent over the field of rational numbers,
then the normal derivative of the harmonic function on SG determined by
$\alpha, \beta, \gamma$ is never zero on any junction point on SG. More
strictly, if there does not exist a relation
$$
n\alpha + m \beta + k \gamma  = 0
$$

\noindent
between $\alpha, \beta, \gamma$ (being not all equal) with $n, m, k$ integers
and $n+m+k=0$, then the normal derivative is never zero on any junction
point on SG.\\

Summarizing the above considerations, we obtain\\

\noindent
{\bf Theorem 5.}\\
{\it The normal derivative of a non-constant harmonic function $f$
on SG is non-zero at all junction points on $G_0 \subset SG$, with at most a
single exception. This exception occurs at a vertex of $G_0$, iff all
restrictions of $f$ to the edges of $G_0$ are monotone.}\\

\centerline{\bf 4. Zero Derivatives}

\noindent
In this section we will show that the derivative of the restriction of a
harmonic function $f$ on SG to an edge of any $G_m$ is differentiable at a
point dividing the edge in ratio 1:3 and the derivative there vanishes.
It is again enough to show this for the edge $[p_1,p_2]=[0,1]$ of $G_0$ as
the extention rule for the harmonic function is the same at every scale.\\

\noindent
{\bf Theorem 6.}\\
{\it Let $f$ be a harmonic function on SG and $p$ the point
dividing the edge $[p_1,p_2]$ in ratio 1:3. (i.e. $p=1/3$). Then
$$
(f|_{[0,1]})^{\prime}(\frac{1}{3})=0.
$$
}
\noindent
{\bf Proof.} Let us first assume that the restriction of $f$ to $[0,1]$
is monotone increasing. To approach the point $p=\frac{1}{3}$ from left
and right with geometrically convergent sequences we use the following
sequence of triangles $\triangle_m=\{ p_0^m, p_1^m, p_2^m \}$:\\

Let $\triangle_0= G_0= \{ p_0, p_1, p_2 \}$ and let $\triangle_m$ be defined
as in Fig 3. (right third of the left third of $\triangle_{m-1}$)\\

\begin{figure}[h]
\begin{center}
\includegraphics{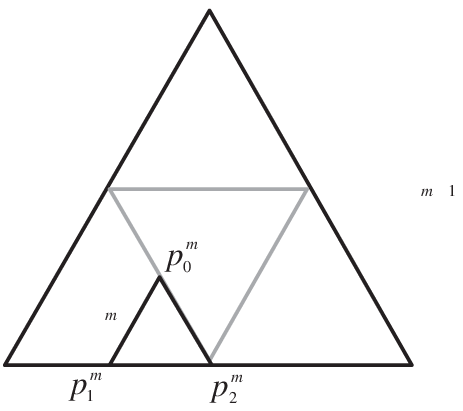}
\caption{The sequence of triangles $\triangle_m$}
\end{center}
\end{figure}

\noindent
One can compute

$$
p_1^m=\frac{1}{4}+(\frac{1}{4})^2+...+(\frac{1}{4})^m=\frac{1}{3}-\frac{1}{3}(\frac{1}{4})^m
$$

$$
p_2^m=p_1^m+(\frac{1}{4})^m=\frac{1}{3}+\frac{2}{3}(\frac{1}{4})^m.
$$

\noindent
Let $f(p_0^m)=\alpha_m, \, f(p_1^m)=\beta_m$ and $f(p_2^m)=\gamma_m, \,
\alpha_0, \beta_0, \gamma_0$ being $\alpha, \beta, \gamma$.\\

\noindent
We want to compute the values $\beta_m$ and $\gamma_m$ explicitly. Using (3)
we get

\begin{equation}
\alpha_m=\frac{1}{25}[6 \alpha_{m-1} + 13 \beta_{m-1} + 6 \gamma_{m-1} ]
\end{equation}

\begin{equation}
\beta_m=\frac{1}{25}[4 \alpha_{m-1} + 16 \beta_{m-1} + 5 \gamma_{m-1} ]
\end{equation}

\begin{equation}
\gamma_m=\frac{1}{5}[\alpha_{m-1} + 2 \beta_{m-1} + 2 \gamma_{m-1} ]
\end{equation}

\noindent
From (13)-(15) we obtain

$$
5\alpha_m+15\beta_m+7\gamma_m = 5\alpha_{m-1}+15\beta_{m-1}+7\gamma_{m-1}
$$

\noindent
for all $m=1,2,...$. In other words,

\begin{equation}
5\alpha_m+15\beta_m+7\gamma_m = 5\alpha+15\beta+7\gamma=:c.
\end{equation}

\noindent
From (16) and continuity of $f$ we get

\begin{equation}
f(\frac{1}{3})=\frac{c}{27}.
\end{equation}

\noindent
Using (16) we can eliminate $\alpha_{m-1}$ from (15):

\begin{equation}
\beta_m=\frac{1}{125}[4c+20\beta_{m-1} - 3 \gamma_{m-1}]
\end{equation}

\begin{equation}
\gamma_m=\frac{1}{125}[5c-25\beta_{m-1} + 15 \gamma_{m-1}]
\end{equation}

\noindent
As can be seen from (18) and (19), the sequence

\begin{equation}
t_m=u\beta_m+v\gamma_m
\end{equation}

\noindent
with $u=10, \quad v=1-\sqrt{13}$, satisfies the recursion formula

\begin{equation}
t_m=w+st_{m-1},
\end{equation}

\noindent
where $w=\frac{9-\sqrt{13}}{25}c, \quad s=\frac{7+\sqrt{13}}{50}$.\\

\noindent
From (21) $t_m$ can be determined:

\begin{equation}
t_m=w\frac{s^m-1}{s-1}+s^m.t_0  \qquad (t_0=10\beta+(1- \sqrt{13})\gamma ).
\end{equation}

\noindent
From (18), (19) and (21) we obtain

$$
\gamma_m=\frac{c}{25}-\frac{1}{50}t_{m-1}+\frac{v+6}{50}\gamma_{m-1} ,
$$

\noindent
and inserting $t_{m-1}$ from (22) we get

\begin{equation}
\gamma_m= l + k.s^{m-1}+h.\gamma_{m-1},
\end{equation}

\noindent
where $l=\frac{c}{25}+\frac{w}{50(s-1)}, \quad
k=-\frac{1}{50}(\frac{w}{s-1}+t_0), \quad h=\frac{v+6}{50}$.\\

\noindent
The recursion (23) gives $\gamma_m$ explicitly:

$$
\gamma_m=[\frac{l}{h-1}-\frac{k}{s-h}+\gamma]h^m+\frac{k}{s-h}s^m+\frac{c}{27}.
$$

\noindent
As $0<h<\frac{1}{4}$ and $0<s<\frac{1}{4}$ we obtain finally

$$
\lim_{m \longrightarrow \infty} \frac{f(p_2^m)-f(\frac{1}{3})}{p_2^m-\frac{1}{3}}=0.
$$

\noindent
Taking $x_0=\frac{1}{3}, \, d=\frac{1}{4}$ and $a=\frac{2}{3}$ in Lemma 3,
we see that the right derivative of the restriction of $f$ to $[p_1,p_2]=[0,1]$
at $p=1/3$ exists and is zero.\\

\noindent
Similarly, from (18), (19), (21) we get

$$
\beta_m=\frac{1}{u}[\frac{w}{s-1}-\frac{vk}{s-h}+t_0].s^m-\frac{v}{u}
(\frac{l}{h-1}-\frac{k}{s-h}+\gamma).h^m+\frac{c}{27}
$$

\noindent
and this shows that the left derivative at $p=\frac{1}{3}$ exists and is also zero.
Together we obtain

$$
(f|_{[0,1]})^{\prime}(\frac{1}{3})=0.
$$

\noindent
Now we consider the case where the restriction of $f$ to $[0,1]$ is not
monotone. In that case we know that the restriction is monotone in two pieces.
If the extremum is not attained at $p=1/3$, then there is a neighborhood
$(\frac{1}{3}-\delta, \frac{1}{3}+\delta )$ where the restriction is
monotone and the above proof applies. If the extremum is attained at $p=1/3$,
then Lemma 2,3 and the above proof works still on two sides of $p=1/3$ and we
get
$
(f|_{[0,1]})^{\prime}(\frac{1}{3})=0 \quad \Box
$\\

Address:\\
Anadolu University\\
Mathematics Department\\
26470 EskiŸehir-TURKEY\\

e-mails:\\
bdemir@anadolu.edu.tr\\
vcaferov@anadolu.edu.tr\\
skocak@anadolu.edu.tr\\
mureyen@anadolu.edu.tr

\end{document}